\documentclass{article}
\usepackage[utf8]{inputenc}
\usepackage{amsmath,amssymb,dsfont,amsfonts}
\usepackage{mathrsfs}
\usepackage{graphics}
\usepackage{graphicx}
\usepackage{color}

\newcommand{\R}{\mathbb{R}}

\newcommand{\diff}{\mathrm{Diff}}
\newcommand{\proof}{{\noindent \bf Proof:} }
\newcommand{\eop }{ \hfill $\Box$ }

\newtheorem{theorem}{Theorem}[section]
\newtheorem{proposition}[theorem]{Proposition}
\newtheorem{corollary}[theorem]{Corollary}

\newtheorem{definition}[theorem]{Definition}
\newtheorem{lemma}[theorem]{Lemma}
\usepackage[bottom=2cm,top=2.5cm,left=2.5cm,right=2.5cm]{geometry}

\title{Geometric decomposition of flows generated by rough path differential equations}
\author{Pedro Catuogno\footnote{Mathematics Department, State University of Campinas, Brazil. E-mail: pedrojc@unicamp.br. 
Research partially supported by São Paulo Research Foundation - FAPESP nr.2020/04426-6.} \ \ \ \ Lourival Lima\footnote{Yachay Tech University, Ecuator. E-mail: lrodrigues@yachaytech.edu.ec. 
Research supported by Yachay Tech University.}\ \ \ \  Paulo Ruffino\footnote{Mathematics Department, State University of Campinas, Brazil. E-mail: ruffino@unicamp.br. 
Research partially supported by Brazilian Research Council - CNPq 305212/2019-2, São Paulo Research Foundation - FAPESP nr.2020/04426-6.}    }

\begin{document}

\maketitle

\begin{abstract} 

Whenever an It\^o-Wentsel type of formula holds for composition of flows of a certain differential dynamics, there exists locally a decomposition of the corresponding flow according to complementary distributions (or foliations, in the case of integrability of these distributions). Many examples have been proved in distinct context of dynamics: Stratonovich stochastic equations, L\'evy driven noise, low regularity $\alpha$-H\"older control functions ($ \alpha\in (1/2,1]$), see e.g. \cite{CLR}, \cite{CSR}, \cite{MMR18}, \cite{MMR16}. Here we present the proof of this categorical property: we illustrate with the $\alpha$-H\"older rough path, $\alpha \in (1/3, 1/2]$ using the It\^o-Wentsel formula in this context proved in \cite{CCM}. Different from the previous approaches, here however, instead of using an intrinsic rough path calculus on manifolds, the manifold has to be embedded in  an Euclidean space.  A cascade decomposition is also shown when we have multiple lower dimensional directions which span the whole space. As application, the linear case is treated in details: the cascade decomposition provides a row factorization of all matrices which allow real logarithm.

\end{abstract}

\bigskip

\bigskip


\vspace{0.3cm}
\noindent {\bf MSC2020 subject classification:} 60L20, 37H05 37C10 (60L90).
\section{Introduction}

This article presents a generalization to dynamical systems generated by rough paths after a series of other results regarding geometric decomposition of flows on manifolds. Previous results have considered dynamics generated by Stratonovich equations \cite{CSR}, \cite{MMR18}, driven by semimartingales with jumps \cite{MMR16} and others \cite{CLR}, all of them was proved using an intrinsic calculus on the manifold. To be more precise, let $\varphi_t$ be a flow of (local) diffeomorphisms generated by an autonomous dynamical system on $M$, a connected differentiable manifold with a prescribed complementary pair of subspaces (distributions) in each tangent space. The key point here is the following idea: if an It\^o-Wentzel (Leibniz) type of formula holds for composition of two families of diffeomosphisms of a certain dynamics (e.g. stochastic, rough path generated or others), then there exists locally a decomposition of the original flow $\varphi_t$ according to the complementary distributions (or foliations, in the case of integrability of these distributions). This geometrical structure will be explained in further details in the sequel. In this context, this article explores the recently proved It\^o-Wentzel formula in \cite{CCM} for $\alpha$-rough path, with $\alpha\in (1/3, 1/2]$. 

\bigskip

Rough path differential equations (RDE, for short) is not only interesting by itself as an important sort of control equation, but it is also relevant since it generalizes results in the classical stochastic case. Recently after many foundational works, see e.g. \cite{Friz and Hairer}, \cite{Gubinelli}, \cite{Bailleul-flows}, \cite{Baudoin} among many others and references therein, there have been an increasing interest of the topic in the literature. Just to  mention few of them closely related to the dynamics we treat here see e.g. \cite{Inahama}, \cite{NNT} with similar geometrical struture we use here.

\bigskip

In the next section we present preliminaries of the rough path calculus for $\alpha\in (1/3, 1/2]$ and the It\^o-Wentzel formula we are going to use in our main decompositon theorem \ref{theorem: main decomposition}. In Section 3 we present the geometrical framework and the proof of the main theorem. Different from previous approach where the geometrical analysis is intrinsic to the manifold, here the manifold is embedded in an Euclidean space, where the rough path calculus holds. Decomposition results are typically local in space and time (see e.g. Example 2 for explosion in time). In section 4 we present corollaries related to applications. In particular, a detailed calculation is provided for the linear case. In the main result of the linear section we show that for any linear flow generated by a RDE, there exists a basis in the Euclidean space $\R^n$ such that the associated flow $\varphi_t$ can be written as a composition of linear functional as:
\[
   \varphi_t = \left( \begin{array}{cccc}
\mathbf{[}* & * \ldots *  & *&   * ]_{d_1\times n} \\
  & & & \\
   & I_{d_2}&  & \\
 &  & \ddots & \\
  &&& \\
 & & & I_{d_k} \\
 \end{array} \right) 
\left( \begin{array}{cccc}
I_{d_1} &  &  &  \\
 & & & \\
 \mbox{{\bf [}} * &  *& * & * \mbox{\bf{]}}_{d_2 \times n} \\
  &  I_{d_3}  & & \\
 &  & \ddots & \\
 & & & I_{d_k} \\
 \end{array} \right) \ldots 
 \left( \begin{array}{cccc}
I_{d_1} &  &  &  \\
&&& \\
  &  \ddots &  &  \\
  &&& \\
  &  & I_{d_{k-1}} & \\
\mbox{\bf{[}} * & * \ldots * & *  & * \mbox{\bf{]}}_{d_k \times n}  \\
 \end{array} \right),
\]
where $d_i=1$ if the corresponding eigenvalues is real, otherwise $d_i=2$, with no explosion in time. The decomposition do not depend on the subspace of generalized eigenvector.


\section{Rough path integration and It\^o-Wentzel Formula}

In this section, we recall basic definitions and properties of rough path theory, essentially only the necessary to develop our theory of decomposition of flows, so that, for details  and further equally interesting properties we recommend, among many others, \cite{Friz and Hairer}, \cite{Gubinelli}, \cite{Baudoin}, \cite{Bailleul-flows}  and references therein.  The crucial property we are going to use in order to obtain the decompostion of rough generated flows is the  It\^o-Wentzel type formula in this context. Another simplification in this rather short  and direct overview on the preliminaries of the theory is the fact that all linear spaces involved here shall be Euclidean spaces. In fact, these finite dimensional spaces are enough to perform our decomposition via embedding differentiable manifolds in a sufficiently large dimensional space.



\bigskip

The interval of time we are considering here is always $[0,T]$, with $0<T<\infty $. The parameter $\alpha$ is always in $(1/3, 1/2]$ and the domain of the tensor product component is the `simplex'  given by  $\Delta_T :=\{ (s,t): 0\leq s \leq t \leq T\}$.

\begin{definition}
An $\alpha$-H\"older rough path  is a pair of functions  $\textbf{X} = (X, \mathbb{X})$, such that $X :[0,T] \rightarrow \R^d$  is $\alpha$-H{\"o}lder continuous and $\mathbb{X}:\Delta_T \rightarrow \R^d \otimes \R^d$ is $2 \alpha$-H{\"o}lder continuous path, such that the following (Chen) relation holds:
\begin{eqnarray}
\mathbb{X}_{st}= \mathbb{X}_{su} + \mathbb{X}_{ut} + X_{su}\otimes X_{ut},
\end{eqnarray}
for all $0 \leq s \leq u \leq t $, where $X_{st}= X_t -X_s$.

\end{definition}

The motivating example for the theory  is lift of the canonical Brownian motion via the L\'evy area, i.e. $\textbf{B} = (B, \mathbb{B})$ where $B_t$ is a  Brownian motion in $\R^d$, and 
\[
 \mathbb{B}_{st} := \int_s^t B_{sr}  \otimes dB_r := \sum_{i,j} \mathbb{B}^{ij}_{st} \ e_i \otimes e_j
\]
where $\{e_ i: 1, \ldots ,  d\} $ is the canonical basis of $\R^d$ with
\begin{equation} \label{Eq: ito-Stratonovich BM}
\mathbb{B}^{ij}_{st} = \int_s^t B_{sr}^i \ dB^j_r.
\end{equation}


We denote the transpose in the tensor product $\R^d \otimes \R^d $ by $(a\otimes b)^* = (b\otimes a)$. Among the space of rough paths with first coordinate $X;[0,T] \rightarrow \R^d$, we are interested in those which satisfy the following:

\begin{definition}
Let $\textbf{X} = (X, \mathbb{X})$, be a rough path. We say that   $\textbf{X}$ is {\it geometric} if: 
\[
 X_{st}\otimes X_{st} = 2\ 
 \mathrm{Sym}(\mathbb{X}_{st}),
\]
where, $\mathrm{Sym}(\mathbb{X}_{st}) := \frac{1}{2}(\mathbb{X}_{st} + \mathbb{X}_{st}^*)$ is the symmetric part of $\mathbb{X}$. 
\end{definition}

For example, the rough path given by the lift of the  Brownian motion $\textbf{B} = (B, \mathbb{B})$ defined above is geometric if the integral of equation (\ref{Eq: ito-Stratonovich BM}) is taken in the sense of Stratonovich, see e.g. \cite{Friz and Hairer}. The importance of the next definition is that {\it controlled paths} are the natural integrands with respect to rough paths:

\begin{definition}
Fix $X: [0,T] \rightarrow \R^d$ an $\alpha$-H\"older continuous path. Let $Y: [0,T] \rightarrow \R^\ell$ and $ Y': [0,T] \rightarrow L(\R^d,\R ^\ell)$ be also $\alpha$-H\"older paths. A pair $(Y, Y')$ is called a {\it  controlled path} with respect to $X$ if 
$$
(s,t) \in \Delta_T \longmapsto R^Y_{st}:= Y_{st} - Y'_s \, X_{st}
$$ 
is $2\alpha$-H\"older continuous. The path $Y'$ is called a {\it Gubinelli derivative} of $Y$ with respect to $X$.
\end{definition}
The prototype of a controlled path with respect to $X$ is given by $Y=F(X)$, where $F: \R^d \rightarrow L(\R^d; \R^m)$ is a differentiable function and $Y' = DF(X)$. 
Next Theorem gives us the definition and the main property of the rough path integral:

\begin{theorem}
 \label{Def: rough integral}
Let  $\textbf{X} = (X, \mathbb{X})$ be an  $\alpha$-H\"older  rough path in $\R^d$ and $(Y, Y' ) $ be a controlled path with respect to $X$, with values in $\R^\ell \cong L(\R^d; \R^m)$. Then the following limit

\begin{eqnarray*}
\int_s^tY_r\ d\mathbf{X}_r := \lim_{|\pi| \to 0}\sum_{[u,v]\in \pi}  \Big( Y_u\, X_{uv} + Y'_u \ \mathbb{X}_{uv} \Big)
\end{eqnarray*}
exists, where $\pi$ are partitions of $[s,t]$ and the limit is taken when the mesh of $ \pi$ goes to zero. 

\end{theorem}
This limit is called the rough integral of the controlled path $(Y, Y')$ with respect to the rough path $\mathbf{X}$. Mind that $Y'$ has values in $L(\R^d; L(\R^d; \R^m)) \cong  L(\R^d \otimes \R^d ; \R^m)$. An important property of this rough integration is the following estimate for all $0\leq s \leq  t \leq T$, 
\[
\Big|\int_s^tY_rdX_r - Y_sX_{st} - Y'_s \mathbb{X}_{st} \Big| \leq C \Big\{ \big| {R}^Y \big|_{2\alpha} \, \big|X \big|_{\alpha} + \big|Y' \big|_{\alpha}\big| \mathbb{X} \big|_{2\alpha} \Big\} |t-s|^{3\alpha},
\]
where $C$ is a constant which depends only on $\alpha$, and  the norms on the right hand side are the standard  $\alpha$-H\"older norms, see e.g. \cite{Friz and Hairer}, \cite{Gubinelli}. The key point in the link between  stochastic analysis and rough path theory is the following: if $(Y, Y' )$ is an  adapted stochastic process which is  $\alpha$-H\"older controlled with respect to $B$  then 
\[
\int_0^T Y_s \, d\mathbf{B}_s = \int_0^T Y_s \, d B_s, 
\]
where the integral on the right hand side is It\^o (Stratonovich) integration if in equation (\ref{Eq: ito-Stratonovich BM}) we have It\^o (Stratonovich, respectively) integral. See e.g. \cite[Prop. 5.1 and Cor. 5.2]{Friz and Hairer}.

\bigskip

The integration of controlled path with respect to another controlled path is a fundamental tool in the rough path calculus. Let $(Y, Y')$ and $(Z,Z')$ be two controlled path with respect to a rough path $X$. Then the limit

\begin{eqnarray*}
\int_s^tY_r\ d Z_r := \lim_{|\pi| \to 0}\sum_{[u,v]\in \pi}  \Big( Y_u\, Z_{uv} + Y'_u \ Z'_u \mathbb{X}_{uv} \Big),
\end{eqnarray*}
exists, where $\pi$ are partitions of $[s,t]$ and the limit is taken when the mesh of $ \pi$ goes to zero. Mind that here we have that $Y$ is a controlled path with values in $L(\R^m, \R^\ell)$ and $Z$ has values in $\R^m$.



\bigskip

In this article we consider the flow of local diffeomorphisms of the dynamics generated by a rough differential equation (RDE). More precisely, let $\mathbf{X}$ be an  $\alpha$-H\"older rough path in $\R^d$ and $F: \R^m \rightarrow L(\R^d; \R^m)$ be a differentiable map (cf. $d$ vector fields in $\R^m$). Consider the RDE given by:

\begin{equation} \label{Eq: dinamica original}
    dx_t = F(x_t) \ d \mathbf{X},
\end{equation}
where the solution is taken in the sense of the integration above. We denote by $\varphi_t$ the solution flow of (local) diffeomorphisms. For existence and uniqueness of solution given and initial conditions and  existence of flow of (local) diffeomorphisms  see e.g. \cite{Bailleul-flows}, \cite[Chap. 8]{Friz and Hairer}. We remark that via the link between stochastic and rough path theory, the calculation of the solutions in both systems are analogous, i.e. if $y_t$ satisfies 
\[
y_t = y_0 + \int_0^t F(y_s) \, d\mathbf{B}_s 
\]
for each $\omega$ (a.s.) then 
\[
y_t = y_0 + \int_0^t F(y_s) \, dB_s, 
\]
where the integration above can be either in It\^o or Stratonovich sense. And reciprocally, see e.g. \cite[Thm 9.1]{Friz and Hairer}.

\bigskip


We finish this section with the generalized It\^o formula and the corollary we are going to use in the next section:

\begin{theorem}[It\^o-Wentzel formula] \label{Theorem: Ito-Ventzel}
Let $\textbf{X} = (X, \mathbb{X})$ be an $\alpha$-H\"older geometric rough path. 

Let 
$x\in \R^\ell \longmapsto (h(\cdot, x), h'(\cdot, x)) $ be a continuous family of controlled rough path with respect to $X$ with values in $L(\R^d,\R^m)$. We assume that   
\begin{enumerate}
\item $h(\cdot, x)$ is continuous and twice differentiable in $x$.
\item $D_xh(\cdot, x)$  is continuous and differentiable in $x$. 
\item For each $t \in [0,T]$, $h'(t,x)$ has continuous derivative in  $x$.
\item $x\mapsto (D_xh(\cdot, x), (D_xh)'(\cdot, x)) $ is a continuous family of controlled rough path with values in $L(\R^\ell; L(\R^d;\R^m))$. 
\end{enumerate}
Let
\[
 g(t,x) = g(0,x) + \int_0^t h(s,x)d\mathbf{X}_s.
\]
Assume that $g: [0,T] \times \R^\ell \rightarrow \R^m$ is continuous and twice differentiable in $\R^\ell$ and the functions  $(t, x) \mapsto  D g(t,x)$  and $(t, x) \mapsto  D^2 g(t,x)$ are continuous.
Then for $(Z, Z')$ a controlled rough path w.r.t. $X$ we have that:
\begin{eqnarray}\label{Formula: IW}
 g(t,Z_t) & = & g(0,Z_0) + \int_0^t h(r,Z_r)d\mathbf{X}_r + \int_0^t Dg(r,Z_r) dZ_r 
\end{eqnarray}
 \end{theorem}

Cf. Kunita's It\^o generalized formula in \cite{Kunita1}, \cite{Kunita2} and \cite{Kunita3}. For a proof in the rough path context, see \cite{CCM}. Two direct corollaries come from this generalized It\^o formula. The first one is the formula for changing  coordinates (classical It\^o formula): 

\begin{corollary}[Change of coordinates] \label{Cor: change of coordinates} For geometric rough path and the same notations of Theorem \ref{Theorem: Ito-Ventzel}. If $h \equiv 0$ and $g: \R^\ell \rightarrow \R^m$ is a differentiable function then
\[
g(Z_t)= g(Z_0) + \int_0^t Dg (Z_r) \, dZ_r
\]
In particular is $Z_t$ is a solution of equation (\ref{Eq: dinamica original}), then 
\[
g(Z_t)= g(Z_0) + \int_0^t Dg F (Z_r) \, d \mathbf{X}_r.
\]
Moreover, if $g$ is a diffeomorphism, then  it is a conjugacy between  $Z_t$ and $g(Z_t)$. 
    
\end{corollary}

\proof It is straightforward.

\bigskip

The second corollary is crucial for the decomposition of flows which appears in the next Section: 

 \begin{corollary} \label{Cor: IW}
Let $Y$ and $Z$ be the flows associated with the RDE  driven by a geometric $\alpha$-H\"older rough path $\mathbf{X}$: 
\[
dY = G(Y) d\mathbf{X}~ \ \ \ \ \mbox{ and }~ \ \ \ \  dZ =H(Z)d\mathbf{X}
\] 
respectively. 
Then the composition of the flows $V= Y \circ Z$ satisfies the RDE:
\[
 dV = G(V)d\mathbf{X} + Y_*H(V)d\mathbf{X}.
\]     
 \end{corollary}
\proof We apply the It\^o-Wentzell formula (\ref{Formula: IW}).
By assumption,
\[
Y_t(x) = x + \int_0^t G(Y_s(x))d\mathbf{X}_s.  
\]
Applying the It\^o-Wentzell formula, the substitution rule and the expression for the Gubinelli derivative of the solution of a rough differential equation we obtain
\begin{eqnarray*}
V_t & = &  x + \int_0^t G(Y_r(Z_r))d\mathbf{X}_s + \int_0^t D_xY_r(Z_r)d Z_r \\
& =& x + \int_0^tG(V_r)d\mathbf{X}_r + \int_0^t D_xY_r(Z_r)H(Z_r)d\mathbf{X}_r \\ 
& = & x + \int_0^t(G(V_r)+ (Y_r)_* H(V_r))d\mathbf{X}_r.
\end{eqnarray*}

\eop

\section{Decomposition of Flows according to complementary distributions}

Let $M$ be a (compact) connected $m$-dimensional  differentiable manifold. Like all the others geometrical objects, differentiability here means smoothness.
To fit the rough path theory presented in the previous sections into rough dynamics in $M$, we shall, firstly, consider $M$ as a somehow nice subset (surface) of an Euclidean space. Indeed, there always exists an immersion $i:M \rightarrow \R^N$  for $N$ large enough: either applying the classical Whitney theorem in general, or Nash isometric theorem, if $M$ is Riemannian. We shall freely identify $M$ with $i(M) \subseteq \R^N$, its image by the immersion. 

\bigskip

A version of the support theorem holds for RDE in the following sense. In equation (\ref{Eq: dinamica original}) if the image $ \mathrm{Im} \{ F(x) \} \subseteq T_xM$ for all $x\in M$, then the solution flow 
$\varphi_t$  in $ \R^N$ preserves $M$, i.e. the solutions starting in $M$ lies in $M$ for all lifetime of the solution, yet:  $\varphi_t (x)\in M$ for all $x\in M$ and $t<T\wedge \tau$, where $\tau$ is the explosion time of the solution. In fact, this is not difficult to prove. Given an initial conditions $x_0 \in M$ let $s_0 := \sup \{ 0\leq s \leq T\wedge \tau : \varphi_t(x_0) \in M \} $. One just have to prove that the interval $[0,s_0]$ is an open set in $[0, T\wedge \tau ]$. Take a local coordinate around $\varphi_{s_0}(x_0)\in M$ such that locally $M$ is mapped horizontally in $\R^m \times \R^{N-m}$. In a neighbouhood of $\varphi_{s_0}(x_0)$, the vertical component in $\R^{N-d}$ of the integrand function $F$ vanishes. Hence, also in a neighborhood of $s_0$, with respect to this coordinates, the solution $\varphi_s(x_0)$ is also horizontal, i.e. $\varphi_s(x_0) \in M$. Hence $s_0 = T\wedge \tau$, i.e. the solution stays always on $M$. 

\bigskip
 Consider a geometric rough path $\mathbf{X}=(X, \mathbb{X})$ in $\R^d$ and let $F:M \rightarrow L(\R^d, TM)$ be a differentiable mapping (`$d$ vector fields') such that $F(x):\R^d \rightarrow T_xM$ for all $x\in M$. A RDE on $M$  written as 
\begin{equation} \label{Eq: dinamica original em M}
    dx_t = F(x_t) \ d \mathbf{X},
\end{equation}
is interpreted as a RDE in the Euclidean space where $M$ is immersed and $F$ extended smoothly in a neighbourhood of $M\subseteq \R^N$. The existence of flow in $\R^N$ and the support theorem in the rough context  allow us to talk about a flow $ \varphi_t$ of local diffeomorphisms in $M$. In the case of a non-autonomous $F(t,x_t)$, the solution flow is denoted by $\varphi_{s,t}$.

\bigskip
Let $\mbox{Diff}(M)$ be the infinite dimensional Lie group of smooth diffeomorphisms of $M$  generated by flows of smooth vector fields on $M$ (understood as extended smoothly vector fields in $\R^N$). The Lie algebra associated to $\mbox{Diff}(M)$ is the infinite dimensional space of smooth vector fields on $M$, see e.g. Neeb \cite{Neeb}, Omori \cite{Omori}, among others. In this Lie group, the exponential map $\exp \{tY\} \in \mbox{Diff}(M)$ is the associated flow of global diffeomorphisms generated by the smooth vector field $Y$. In this context, given an element $\varphi\in \mbox{Diff}(M)$ the derivative of the right translation is given by $R_{\varphi *}Y= Y(\varphi)$ for any smooth vector $Y$; the derivative of left translation $L_{\varphi *}Y= D\varphi (Y)$, and the adjoint operator is given by $\mathrm{Ad}(\varphi)Y = \varphi_* (Y (\varphi^{-1}))$.

\bigskip

 This Lie group structure is convenient for obtaining the equations of each component of the decomposition we are interested in. In particular, a solution flow $\varphi_t$ of the RDE (\ref{Eq: dinamica original em M}) can be written as the solution of a right invariant RDE in the Lie group of diffeomorphisms $\mbox{Diff}(M)$: 
\begin{equation}
    \label{Eq: right inv}
d\varphi_t =  R_{\varphi_{t*}}Y\, d\textbf{X}_r.
\end{equation}


Interesting problems arise when one decomposes a (flow of) diffeomorphism  $\varphi \in \mbox{Diff}(M)$, into composition of conveniently prescribed components, in the sense that each component provides a geometrical or dynamical information. This kind of decomposition appears in the literature, for example, in Bismut \cite{Bismut}, Kunita \cite{Kunita1} and many others. In particular, it is also relevant when each component of the decomposition belongs to prescribed subgroups of $\mbox{Diff}(M)$, see e.g Melo et al \cite{MMR18}, Catuogno et al \cite{CSR}, Iwasawa and non-linear Iwasawa decomposition \cite{Colonius and Ruffino}, Ming Liao \cite{Liao} among many others.

\bigskip

In this section, we explore in the context of rough path calculus the existence of a geometrical decomposition as treated in \cite{CSR}, \cite{MMR18}, \cite{CLR}, i.e. in the following sense:  Suppose that locally $M$ is endowed with a pair of  differentiable distributions:  i.e., every point $x\in M$ has a neighbourhood $U$ and  differentiable mappings $\Delta^1: U \rightarrow Gr_{k}(M)$ and $\Delta^2: U \rightarrow Gr_{{m-k}}(M)$ respectively, where 
$$
Gr_p (M) = \bigcup_{x \in M}Gr_p(T_x M)
$$ 
is the Grasmannian bundle of $p$-dimensional subspaces over $M$, with $1\leq p\leq m$. We assume that $\Delta^1$ and $\Delta^2$ are complementary in the sense that $\displaystyle{\Delta^1(x) \oplus \Delta^2(x) = T_xM}$, for all $x \in U$.   With this notation we define the subgroup of $\mbox{Diff}(M)$ which is generated by a certain distribution $\Delta$ by:


\begin{equation*}
\mbox{Diff}(\Delta, M) = \mbox{cl} \bigl\{ \exp (t_1X_1) \ldots \exp (t_nX_n), \mbox{ with } X_i\in\Delta, t_i\in\mathbb{R}, \forall n\in\mathbb{N} \bigr\}.
\end{equation*}

 Note that if a distribution $\Delta$ is involutive, then each element of the group $\mbox{Diff}(\Delta, M)$ preserves the leaves of the corresponding foliation.  

\bigskip


\noindent {\bf Extending the distributions to $\R^N$.} We freely treat these complementary distribution on $M$ as intrinsic to $M$ since the complementary distributions $\Delta^1 \oplus \Delta^2 = T_xM$ has locally differentiable  extensions to neighbourhoods of $M$ in $\R^N$: in fact, if $v_1, \ldots, v_k$ are l.i smooth vector fields which span $\Delta^1$, then, there exists locally a neighborhood in $\R^N$ where their extensions span to an extended distribution of $\Delta^1$ to $\R^N$. Here, we are abusing notation when writing the distributions on $M$ and in $\R ^N$ with the same symbol. With this construction, locally the extended $\Delta^1$ and $\Delta^2$ in $\R^N$ will also be transversal. One only has to increase the dimension of $\Delta^2$, say, by the span of vector fields transversal to $M$ to get completeness in $\R^N$, i.e. $\Delta^1 \oplus \Delta^2 = \R^N$. When necessary to  distinguish the $(N-k)$-dimensional $\Delta^2$ extended distribution in $\R^N$ and the original $(m-k)$-dimensional $\Delta^2$ which is only defined on $M$, the last one will be denoted by $\Delta^2_M$.

\bigskip

In particular, in this Section we focus on the subgroups $\mbox{Diff}(\Delta^1, M)$ and $ \mbox{Diff}(\Delta^2, M)$. The main result of  this paper (Theorem \ref{theorem: main decomposition}) establishes a local decomposition of the solution flow $\varphi_t$  of equation (\ref{Eq: dinamica original em M}) into two components: a curve (solution of an autonomous RDE) in $\mbox{Diff}(\Delta^1, M)$ composed with a solutions of a non-autonomous RDE in $\mbox{Diff}(\Delta^2, M)$. 

\begin{definition}
\label{definition: transversality}
A diffeomorphism $\eta \in \mbox{Diff}( M)$ preserves transversality of  $\Delta^1$ and $\Delta^2$ in a neighbourhood $U \subset M$ if  $\eta_{*} \Delta^2 \left(\eta^{-1}(x)\right)\cap \Delta^1(x) = \{0\}$, for all $x \in U$.
\end{definition}

By continuity, given any pair of complementary distributions and any point $p\in M$, there always exists a neighbourhood of the identity $Id \in \mbox{Diff}(M)$ and a neighbourhood $U$ of $p$ where all diffeomorphisms in this neighbourhood preserve transversality. Moreover, if the distribution $\Delta^1$ is involutive then all elements in $\mbox{Diff}(\Delta^1, M)$ preserves tranversality of $\Delta^1$ and $\Delta^2$: in fact, the derivative $\eta_*$ above is a linear isomorphism which sends tangent spaces of the associated foliation to tangent spaces of the same leaf. 

 Next theorem states the main result of this section:

\begin{theorem}[Decomposition of flows of RDE]
\label{theorem: main decomposition}
Locally, up to an explosion time $\tau \in [0,T]$, the solution flow $\varphi_t$ of the RDE (\ref{Eq: dinamica original em M}) is decomposable as
\[
 \varphi_t = \eta_t \circ \psi_t,
\]
where $\eta_t$ is solution of an (autonomous) RDE  in the Lie group $\emph{Diff}(\Delta^1, M)$ and $\psi_t$ is
solution of a non-autonomous RDE in $\emph{Diff}(\Delta^2, M)$.
\end{theorem}

\proof

Given $p \in M$, initially extend the distributions $\Delta^1$ and $\Delta^2$ to a neighbourhood $U$ of $p$ in $\R^N$ as described above.
Take $\eta \in \mbox{Diff}(\Delta^1,\R^N)$ sufficiently close to the identity such that $\eta$ preserves tranversality in $U$, i.e. $\mathrm{Ad}(\eta)\Delta^2$ and $\Delta^1$ are complementary. Mind that $\Delta^1$ might not be involutive, but at points on $M$, the Lie bracket $[\Delta^1, \Delta^1] \subseteq \Delta^1 \oplus \Delta^2$. It means that diffeomrophisms in $ \mbox{Diff}(\Delta^1,\R^N) $ restricts to (local) diffeomorphisms on $M$, i.e. in another notation $\mbox{Diff}(\Delta^1,\R^N)|_M = \mbox{Diff}(\Delta^1,M)$.

\bigskip

The integrand mapping $F$ of the RDE (\ref{Eq: dinamica original em M}), after being extended to $\R^N$,  can be decomposed uniquely as 
\begin{eqnarray} \label{Eq: h and V}
F(x) = h(x)+V(\eta,x),
\end{eqnarray}
where $h(x) \in \Delta^1(x)$ and $V(\eta, x) \in \mathrm{Ad}(\eta)\Delta^2(x)$, for all $x \in U$. Note that for points $x\in M$, $V(\eta,x)\in \Delta^2|_M$, since $F$ itself lives in $\Delta^1 \oplus \Delta^2|_M$ and $\eta_*$ preserves the tangent bundle $TM$. The first component of the decomposition in the statement  of the theorem, $\eta_t$ is taken as the solution of the following RDE in the group $\mathrm{Diff}(\Delta^1, M)$:

\begin{eqnarray}
\label{Eq: para eta geral}
d\eta_t = R_{\eta_{t*}}\, h\, d \mathbf{X}_t,
\end{eqnarray}
with initial condition $\eta_0=Id$, the identity.
Even though the equation above is described in terms of a right translation, it is not a right invariant equation since  $h$ in general depends on $\eta_t$. The second component of the  decomposition of $\varphi_t$ is obtained using that  $\psi_t = \eta^{-1}_t \circ \varphi_t$. Since the families of diffeomorphisms $\varphi_t$ and $\eta_t$ preserves $M$ (i.e. are diffeomorphisms in $M$ when restrict to this space), then $\psi_t$ also preserves $M$.  Applying Corollary \ref{Cor: IW}, it follows that:
\begin{eqnarray*}
d\eta^{-1}_t = -L_{\eta^{-1}_{t*}}\, h\, d\mathbf{X}_t,
\end{eqnarray*}
where $L_{\eta^{-1}_{t*}}$ is the derivative of the left translation at the identity by $\eta^{-1}_t$. Finally, we find a equation for $\psi_t$ applying Corollary \ref{Cor: IW} again:
\begin{eqnarray} 
d\psi_t &=& (\eta^{-1}_t h \, \eta_t \, \psi_t - \eta^{-1}_t X\, \eta_t \, \psi_t)\, d\mathbf{X}_t \nonumber\\
&& \nonumber \\
&=& \mathrm{Ad}(\eta^{-1}_t)V(\eta_t)\ d\mathbf{X}_t. \label{Eq: para psi geral}
\end{eqnarray}
Note that $V(\eta, \cdot)$ does not necessarily belong to $\Delta^2$, but rather  $d\psi_t \in \Delta^2$ since $d\psi_t \in \mathrm{Ad}(\eta^{-1})\mathrm{Ad}(\eta)\Delta^2 = \Delta^2$. Then $\psi_t$ is the $\Delta^2$-component of the original flow $\varphi_t$ in $\R^N$. The intrinsic result for flows on $M$ itself is obtained by restricting the (flow of) diffeomorphisms to $M$, since both $\eta_t$ and $\psi_t$ preserve $M$.

\eop

\begin{corollary}
If the distributions $\Delta^1$ and $\Delta^2$  are integrable, then the decomposition of Theorem \ref{theorem: main decomposition} is unique. 
\end{corollary}
\proof
In fact, in this case $\mathrm{Diff}(\Delta^1, M) \cap \mathrm{Diff}(\Delta^2, M)= \{Id\}$.

\eop

\section{Cascade decomposition and the linear case}

When the state space $M$ is endowed with many compatible pairs of complementary distributions in its tangent bundle (in a sense which will be described below) Theorem \ref{theorem: main decomposition} can be applied recursively to obtain further geometrical applications of the decomposition. In this Section we deal with this case to obtain a sequence of dynamically meaningful  corollaries. In particular, we calculate explicitly the decomposition established in these corollaries for the linear RDE. We prove in particular that the main geometrical obstruction for the decomposition, namely the explosion time, does not appear in this case if one chooses appropriate linear (affine) foliations.    
\bigskip

We start with the assumption that the base manifold $M$ has more than one pair of complementary distributions or foliations. More than that, initially, we assume that there exists  a pair of {\it complementary flag bundles} in the following sense:  There exists a sequence of complementary non trivial distributions, i.e. for all $x\in M$ we have  $\Delta^1_1 (x)\oplus \Delta^2_1(x)= \Delta^1_2\oplus \Delta^2_2= \ldots = \Delta^1_k\oplus \Delta^2_k = T_xM$, such that the first component is an increasing sequence of tangent subspaces $\Delta^1_1 \subset \Delta^1_2 \ldots  \subset \Delta_k^1 \subset T_xM$  and the second component is a decreasing sequence of subspaces  $T_xM \supset \Delta^2_1 \supset \Delta^2_2 \ldots \supset  \Delta_k^2$.

\bigskip

 This structure of nested distribution $F^1 = (\Delta^1_1, \ldots , \Delta^1_k)$ and $F^2=(\Delta^2_k, \ldots , \Delta^2_1)$ are called  {\it flag bundle} over the manifold $M$. Each fibre is a flag manifold, i.e. a homogeneous space $Sl(n, \R)/P$, where $P$ is a parabolic subgroup, see e.g. the classical book \cite{helgason} or more recent applications in dynamics in the survey \cite{San Martin}. If $\dim \Delta^1_i - \dim \Delta^1_{i-1} =1$ for all $i=2, \ldots, k$ then each fibre is called a {\it maximal flag manifold}, in this case $k=m-1$, where $m$ is the dimension of the manifold $M$. Our first corollary establishes a cascade of decomposition in this geometrical context. A version of this result  has already appeared in \cite{CSR} in the context of stochastic systems. We include it here again, in rough path context with a simpler proof.

\bigskip

\begin{corollary} \label{Cor: decomp flag} Given complementary flag bundles $F^1 = (\Delta^1_1, \ldots , \Delta^1_k)$ and $F^2=(\Delta^2_k, \ldots , \Delta^2_1)$  there exists (locally) a decomposition of the original flow, up to an explosion time:

\begin{equation}
    \label{Eq: decomp para flag}
\varphi_t = \eta^1_t \circ \ldots \eta^i_t \circ \eta^{i+1}_t \circ \ldots \eta^k_t \circ \psi_t,
\end{equation}
with the property that for all $i=1, \ldots k$ we have that the first $i$-th components
$(\eta^1_t \circ \ldots \eta^i_t)$ solves a rough differential equation in the Lie group $  \diff(\Delta^1_i)$ and the last $(k-i+1)$ components $(\eta^{i+1}_t \circ \ldots \eta^k_t \circ \psi_t)$ solves a (non autonomous) rough equation in the Lie group $  \diff(\Delta^2_{i})$, with the convention that if $i=k$ then  $\eta^{i+1}_t=Id$. 

The decomposition is unique if the distributions of the flag bundles are all involutive.
\end{corollary}

\proof

By Theorem \ref{theorem: main decomposition}, for each pair of complementary distributions $\Delta^1_i \oplus \Delta^2_i$ there exists a corresponding decomposition 
\[
\varphi_t = \widetilde{\eta_t}^i \circ \widetilde{\Psi_t}^i,
\]
with $\widetilde{\eta_t}^i \in \diff (\Delta^1_i)$ and $\widetilde{\psi_t}^i \in \diff (\Delta^2_{i})$, for  $i=1, \ldots k$. 

\bigskip

For the first and the last elements in the right hand side of formula (\ref{Eq: decomp para flag}) of the statement take  $\eta_t^1 = \widetilde{\eta_t}^1$ and $\psi_t = \widetilde{\psi_t}^k$. And for  $i=2, \ldots k$, define
\[
\eta^i_t := (\widetilde{\eta_t}^{i-1})^{-1} \circ \widetilde{\eta_t}^i.
\]
The result follows immediately by the properties of the decomposition established by Theorem \ref{theorem: main decomposition} and  telescope cancellation.

\eop

\bigskip

A particular situation of complementary flag bundles occurs if the subspaces of the flags are both prescribed by a common sequence of distributions. Namely, suppose we have a sequence of distribution $\Delta_1, \ldots \Delta_k$ such that for each $x\in M$ we have
that $\Delta_1 \oplus \ldots \Delta_k=T_xM$. 
In this case, there exists a special choice of complementary flag bundles $F_1$ and $F_2$ such that $F_1$ is given by subspaces $\Delta^1_i = \Delta_1\oplus \ldots \Delta_i$ and $F_2$ is constructed with the subspaces $\Delta^2_{i} = \Delta_{i+1}\oplus \ldots \Delta_k$, for $i=1, \ldots, k-1$.. In this case the last component $\psi_t$ in equation (\ref{Eq: decomp para flag}) disappears and the decomposition assumes the following form:


\begin{corollary} \label{Cor: decomp distrib cascata} Consider a sequence of $k$ complementary nondegenerate distributions in a manifold $M$, i.e. $\Delta_1 \oplus \ldots \Delta_k=T_xM$ for all $x\in M$. Then a flow generated by the RDE (\ref{Eq: dinamica original em M}) can be locally decomposed, up to an explosion time as:

\[
\varphi_t = \xi^1_t \circ \ldots \xi^i_t\circ  \xi^{i+1}_t \circ \ldots \xi^k_t,
\]
such that for each $i=1, \ldots k$ we have that the first $i$-th components
$(\xi^1_t \circ \ldots \xi^i_t)$ solves a rough differential equation in the Lie group $  \displaystyle{ \diff(\oplus_{j=1}^i \Delta_j)}$ and the last $(k-i)$ components $(\xi^{i+1}_t \circ \ldots \xi^k_t )$ solves a (non autonomous) rough equation in the Lie group $  \displaystyle{ \diff(\oplus_{j=i+1}^k \Delta_j)}$.

\end{corollary}

\proof It follows directly from Corollary \ref{Cor: decomp flag}: just define the flags $F^1$ and $F^2$ by the following increasing/decreasing non trivial subspaces:
\[
   \Delta^1_i = \bigoplus_{j=1}^i \Delta_j
\]
and 
\[
\Delta^2_i = \bigoplus_{j=i+1}^k \Delta_j.
\]
Mind that, here, the flags $F^1$ and $F^2$ 
are, each one, determined by $(k-1)$ non trivial subspaces.  

\eop

\bigskip

In general we have that the subgroups of diffeomorphisms  $\diff (\Delta^1_i) \subseteq \diff (\Delta^1_{i+1})$ and $\diff (\Delta^2_{i+1}) \supseteq \diff (\Delta^2_{i})$ for all $i=1, \ldots, (k-1)$.  In the case of flags with only involutive distributions, these inclusions are strict, which allows a more refined result:

\begin{corollary} \label{Cor: decomp distrib cascata involutiva} Consider a sequence of $k$ complementary involutive distributions (foliations) in a manifold $M$, i.e. the tangent subspace of the $k$ foliations  $\Delta_1 \oplus \ldots \Delta_k=T_xM$ for all $x\in M$. Then the flow generated by the RDE (\ref{Eq: dinamica original em M}) can be locally decomposed, up to an explosion time, as:

\begin{equation} \label{Eq: decomp distr involutiva}
    \varphi_t = \xi^1_t \circ \ldots  \xi^k_t,
\end{equation}
such that for each $i=1, \ldots k$ we have that
$\xi^i_t \in \diff (\Delta_i)$. The decomposition is unique. 
\end{corollary}

\proof

Consider the same flags generated by subspaces  $\Delta^1_i$ and $\Delta^2_i$ as in Corollary \ref{Cor: decomp distrib cascata} above. Again, by Theorem \ref{theorem: main decomposition}, for each pair of complementary distributions $\Delta^1_i \oplus \Delta^2_i$, with $i=1, \ldots, k$, there exists a corresponding decomposition of the same original flow of equation (\ref{Eq: dinamica original em M}), that is: 
\begin{equation}
    \label{eq: pre distr involutive}
\varphi_t = \widetilde{\eta_t}^1 \circ \widetilde{\Psi_t}^1 =\widetilde{\eta_t}^2 \circ \widetilde{\Psi_t}^2 = \ldots \widetilde{\eta_t}^k \circ \widetilde{\Psi_t}^k.
\end{equation}
with $\widetilde{\eta_t}^i \in \diff (\Delta^1_i)$ and $\widetilde{\psi_t}^i \in \diff (\Delta^2_{i})$, for  $i=1, \ldots k$. 

\bigskip

The first and the last components in the equation (\ref{Eq: decomp distr involutiva}) of the statement are taken as  $\xi_t^1 = \widetilde{\eta_t}^1$ and $\xi^k_t = \widetilde{\psi_t}^k$. Obviously $\xi^1_t \in \diff (\Delta_1)$ and $\xi^k_t \in \diff (\Delta_k)$. For  $i=2, \ldots k$, define
\[
\xi^i_t := (\widetilde{\eta_t}^{i-1})^{-1} \circ \widetilde{\eta_t}^i \in \diff (\Delta^1_i). 
\]
Equation (\ref{eq: pre distr involutive}) says that it is also equivalent to define $\xi^i_t=(\widetilde{\Psi_t}^{i})^{-1} \circ \widetilde{\Psi_t}^{i-1} \in \diff (\Delta^2_{i-1})$. One concludes that $\xi^i_t \in  \diff (\Delta^1_i) \cap \diff (\Delta^2_{i-1})$. Since all the distributions are associated to foliations, diffeomorphisms in the intersection $\diff (\Delta^1_i) \cap \diff (\Delta^2_{i-1})$ can only be generated by flows associated to vector fields on the foliation associated to the involutive distribution $\Delta^1_i \cap \Delta^2_{i-1} = \Delta_i$, which implies that $\xi_t^i \in \diff (\Delta_i)$.

\eop

\bigskip

\noindent {\bf Example 1:} In general, let  $\phi=(\phi_1, \ldots,
\phi_m): U \subset M \rightarrow \R^m$ be a local coordinate of $M$. Then $\phi$ establishes in $U$, an open set, $m$ 1-dimensional complementary foliations and corresponding distributions which fulfill the geometrical hypothesis of the previous corollaries. Hence, given the  flow of (local)
diffeomorphisms $\varphi_t$, up to an explosion time, we can decompose locally
\[
 \varphi_t=\xi^1_t \circ  \xi_t^2\ldots \circ \xi_t^m
\]
where each diffeomorphism $\xi^i_t$
preserves the $j$-th coordinates for all $j\neq i$.  The 
composition $(\xi^1_t \circ \xi_t^2 \ldots \circ \xi^i)$
is a solution of a rough differential equation in the groups of diffeomorphisms which preserves
coordinates $\phi_{i+1}, \phi_{i+2}, \ldots , \phi_m$. On the other hand,
the composition $(\xi^{i+1}_t \circ \ldots  \circ \xi^{m})$
preserves coordinates $\phi_{1}, \phi_{2}, \ldots , \phi_i$. The decomposition is unique.

\eop

\bigskip

\subsection{Linear case} 

Consider the case of an autonomous linear RDE in $\R^m$:
\begin{equation} \label{eq: linear}
 d \, x(t) = A\, x(t) \ d \mathbf{X}_t,
\end{equation}
where $A$ is an $m\times m$ real matrix. Let $\varphi_t$ be the associated linear flow in $\R^m$. As for the geometrical structure, let $F_1 \oplus F_2 = \R^m$ be two complementary vector subspaces of $\R^m$. By translations, they generate the complementary affine foliations 
$$
\displaystyle{\mathcal{F}_ 1 = \bigcup_{x\in \R^m} F_1 + x} \ \ \ \ \ \mbox{ and } \ \ \ \ \ \ \mathcal{F}_2 = \bigcup_{x\in \R^m} F_2 + x.
$$
In this Section these foliations are the focus of our attention. For the Cartesian decomposition  $\R^m = \R ^k \times \R^{\ell}$, with $\ell = (m-k)$, the corresponding affine foliations generated by $\R^k$ and $\R^\ell$ are called the the {\it Cartesian foliation}. The concern about pair of affine foliations for linear systems comes from the following:

\begin{lemma}
\label{Lemma: linear com folheacao afim}
 Complementary foliations  $\mathcal{F}_1$ and $\mathcal{F}_2$ in $\R^m$ are affine if and only if for all linear flow $\varphi_t$, the local decomposition of Theorem \ref{theorem: main decomposition} 
\[
\varphi_t = \eta_t \circ \psi_t 
\]
is such that the components $\eta_t$ and $\psi_t$ are also linear. In this case, the domain of the decomposition is $\R^m$.
\end{lemma}

\proof
Since any pair of complementary affine foliation can be obtained by changing bases of the Euclidean space, without lost of generality we can prove only the case of Cartesian foliations $\R^k\times \R^\ell$ in $\R^m$. We have to prove that $\eta_t$ and $\psi_t$ are linear. In fact, in this case, we write the linear solution flow as 
\[
\varphi_t = \left(  \begin{array}{ll}
\Big( F_1 (t) \Big)_{k\times k}     & \Big( F_2(t)\Big)_{k \times \ell} \\
    & \\
    \Big( F_3 (t)\Big)_{\ell \times k} & \Big( F_4(t)\Big)_{\ell \times \ell}
\end{array}
\right).
\]
 Since $\eta_t$ can not change the last $\ell$ coordinates, and $\psi_t$ can not change the first $k$ coordinates, necessarily $\psi_t$ is linear and invertible (for $t$ close enough to zero) of the form:
 \begin{equation} \label{eq: psi component}
 \psi_t = \left(  \begin{array}{ll}
\Big( Id \Big)_{k\times k}     & \Big( 0 \Big)_{k \times \ell} \\
    & \\
    \Big( F_3 (t)\Big)_{\ell \times k} & \Big( F_4(t)\Big)_{\ell \times \ell}
\end{array}
\right).
\end{equation}
Which implies that $\eta_t = \varphi_t \circ \psi_t^{-1}$ is also linear. Conversely, suppose that $\eta_t$ and $\psi_t$ are both linear. Then  $\psi_t$ must have the form  (\ref{eq: psi component}) and similarly, the first component must have the form: 
 \begin{equation} \label{eq: eta component}
 \eta_t = \left(  \begin{array}{ll}
\Big( G_1 (t) \Big)_{k\times k}     & \Big( G_2 (t) \Big)_{k \times \ell} \\
    & \\
    \Big( 0 \Big)_{\ell \times k} & \Big( Id \Big)_{\ell \times \ell}
\end{array}
\right).
\end{equation}
We have to prove that $\mathcal{F}_1$ and $\mathcal{F}_2$ are affine foliations. In fact, for each point $(x_1, x_2)\in \R^k \times \R^\ell$, we prove that the leaves of $\mathcal{F}_1$ and $\mathcal{F}_2$ are affine subspaces passing through $(x_1, x_2)$, each of them parallel to $\R^k \times \{0\}$ and $\{0\}\times \R^\ell$ respectively. Assume initially that $x_2\in \R^\ell$ is non zero. Write the matrix of coefficients as 
\begin{equation} \label{eq: decomp matriz A}
A = \left(  \begin{array}{ll}
\Big( A_1 \Big)_{k\times k}     & \Big( A_2 \Big) \\
    \Big( A_3 \Big) & \Big( A_4 \Big)_{\ell \times \ell}
\end{array}
\right)
\end{equation}
Take zero entries for $A_1, A_2$ and $A_3$. For an $\alpha$-H\"older rough path $X:[0,T] \rightarrow \R$, the action of the second component is reduced to $\psi_t (x_1, x_2)= (x_1, \exp (X_t\, A_4) x_2) $ which, by definition, carries  $(x_1, x_2)$ along the foliation $\mathcal{F}_2$. Since $X_t$, $A_4$ are arbitrary and $x_2\neq 0$ then $\exp (X_t\, A_4) x_2$ can reach any point in $\R^\ell$ at finite time. This proves that the leaves of  $\mathcal{F}_2$ passing through any $(x_1, x_2)$ with $x_2\neq 0$ is the  affine subspace generated by $\{0\}\times \R^\ell$. The similar argument shows that the leaves of $\mathcal{F}_1$ passing through any $(x_1, x_2)$ with $x_1\neq 0$ is the  affine subspace generated by $\R^k\times \{0\}$.

There are several ways to conclude the statement for points $(x_1, x_2)$ with $x_1=0$ or $x_2=0$. One of the simplest is to use a local foliated coordinate system at each of these points and use continuity of the leaves.

\eop 

\bigskip

\noindent {\bf Example 2:} (Linear case) 
Lemma \ref{Lemma: linear com folheacao afim} above, together with Corollary \ref{Cor: decomp distrib cascata involutiva} say that Example 1 above has the following form in the linear case when one consider the Cartesian foliations: up to  a lifetime we have that 

\begin{equation} \label{Eq: linear c explosao}
   \varphi_t = \left( \begin{array}{cccc}
[* & * \ldots *  & *&   * ]_{1\times n} \\
  & & & \\
   & 1 &  & \\
 &  & \ddots & \\
  &&& \\
 & & & 1
 \end{array} \right) 
\left( \begin{array}{cccc}
1 &  &  &  \\
 & & & \\
 \mbox{[} * & * \ldots *  & *&   * \mbox{]}_{1\times n}\\
  &   1& & \\
 &  & \ddots &  \\
 & & & \hspace{0mm}1
 \end{array} \right) \ldots 
 \left( \begin{array}{cccc}
1 &  &  &  \\
&&& \\
  &  \ddots &  &  \\
  &&& \\
  &  & 1 & \\
\mbox{[} * & * \ldots * & *  & * \mbox{]}_{1\times n} 
 \end{array} \right).
\end{equation}

\bigskip

\noindent 
A simple example which illustrates the explosion is the pure rotation on $\R^2$ with the Cartesian foliations: Let $\mathbf{X}$ be a one dimensional rough path and consider the linear system:
\[
dx_t = \left( \begin{array}{cc} 0 & -1 \\ 1 & 0 \end{array} \right) x_t\, d\mathbf{X}_t,
\]
where the flow decomposes as:

\bigskip

\begin{eqnarray}
\label{Ex: explosion}
\hspace{-1cm}\left( \begin{array}{cc} \cos X_t & -\sin X_t \\ \sin X_t & \cos X_t \end{array} \right) &=& \left( \begin{array}{cc} \sec X_t & -\tan X_t \\ 0 & 1 \end{array} \right)\left( \begin{array}{cc} 1 & 0 \\ \sin X_t & \cos X_t \end{array} \right). 
\end{eqnarray}

\bigskip

\noindent Explosion in the coefficients of the decomposition occurs if 
$X_{t} \in \{\frac{\pi}{2}, -\frac{\pi}{2}\}$. 

\eop

\bigskip

Next proposition shows that the explosion of Example 2 does not appear in generalized eigenvectors related to real eigenvalues. In such a way that in higher dimensions, isolating subspaces of dimension at most 2 (for nonreal eigenvalues), we have decomposition without explosion. In other words, there always  exists a basis in $\R^n$ such that the decomposition of  equation (\ref{Eq: linear c explosao}) holds for all time $t\in [0,T]$ if we allow 2-dimensional foliations related to non real eigenvalues. More precisely:

\begin{proposition} \label{Prop: decomp no explosion} Given a canonical Jordan basis of the matrix of coefficients (with 1 above the diagonal), the following cascade decomposition of the linear flow $\varphi_t$ has no explosion with respect to this basis:
  \[
   \varphi_t = \left( \begin{array}{cccc}
\mathbf{[}* & * \ldots *  & *&   * ]_{d_1\times n} \\
  & & & \\
   & I_{d_2}&  & \\
 &  & \ddots & \\
  &&& \\
 & & & I_{d_k} \\
 \end{array} \right) 
\left( \begin{array}{cccc}
I_{d_1} &  &  &  \\
 & & & \\
 \mbox{{\bf [}} * &  *& * & * \mbox{\bf{]}}_{d_2 \times n} \\
  &  I_{d_3}  & & \\
 &  & \ddots & \\
 & & & I_{d_k} \\
 \end{array} \right) \ldots 
 \left( \begin{array}{cccc}
I_{d_1} &  &  &  \\
&&& \\
  &  \ddots &  &  \\
  &&& \\
  &  & I_{d_{k-1}} & \\
\mbox{\bf{[}} * & * \ldots * & *  & * \mbox{\bf{]}}_{d_k \times n}  \\
 \end{array} \right),
\]
where $d_i=1$ if the corresponding eigenvalues is real, $d_i=2$ otherwise.
\end{proposition}
   
The positive integer $k$ above counts the  number of real eigenvalues plus the number of pairs of conjugate eigenvalues counting multiplicity, such that $\sum_{i=1}^k d_i= m$.
  
 \bigskip
  
   \proof
We use the notations established  in equations (\ref{eq: psi component}), (\ref{eq: eta component}) and (\ref{eq: decomp matriz A}) for the components $\psi_t$, $\eta_t$ and the matrix of coefficients $A$, with respect to a single Cartesian decomposition. The differential equations of $\psi_t$ and $\eta_t$ are given in general by equations (\ref{Eq: para eta geral}) and (\ref{Eq: para psi geral}) in the proof of the decomposition theorem \ref{theorem: main decomposition}. In the particular case of a linear vector field and Cartesian foliations we have that the submatrices of $\eta_t$ and $\psi_t$ satisfies:

\begin{eqnarray}
\begin{array}{ll}
d G_1(t) = \big[  A_1 \ G_1(t) - G_2(t)\ A_3 \ G_1(t) \big] \, dX_t,    \ \ \ \ \ \ \ \  & \\
 & \\
d G_2(t) = \big[ A_1 G_2(t)+ A_2 - G_2(t)\, A_4 - G_2(t)\, A_3\, G_2(t) \big] \, dX_t,  \ \  & \\
  & \\
d F_3(t) = \big[ A_3\,G_1(t) + A_3\, G_2(t)\, F_3(t) + A_4 F_3(t) \big] \, dX_t,     &   \\
  & \\
  d F_4(t) = \big[ A_3\, G_2(t)\, F_4(t) + A_4\, F_4(t) \big] \, dX_t. & 
\end{array} 
\label{Eq: Eq: constituent g_1 g_2 f_3 f_4}
\end{eqnarray}
For details of this calculation see \cite[Eq. (22)]{CLR}. It is well known that $F_3(t)$ and $F_4(t)$ have no explosion. Potential explosion happens in the equations of $G_1(t)$ and $G_2(t)$. Moreover, an explosion may happen due to the second order terms which here are all associated to the submatrix $A_3$. Hence, there is no explosion if the lower left submatrix $A_3= 0$. In Example 2 above note that explosion happens due to $A_3=[1]$.  

\bigskip

The canonical Jordan form of matrix $A$ give us exactly $k-1$ distinct decomposition of $\R^m = \R^s \times \R^\ell $ with $s= d_1+ \ldots d_j$ and $\ell = d_{j+1} + \ldots d_k$
such that the $A_3$ submatrices of $A$ with respect to this basis vanishes. Corollary \ref{Cor: decomp distrib cascata involutiva} says that there exists the corresponding $k$-factor decomposition. Lemma \ref{Lemma: linear com folheacao afim} guarantees that all the components are linear. And finally equations (\ref{Eq: Eq: constituent g_1 g_2 f_3 f_4}) guarantees that none of the $k$ factors of decomposition explodes at finite time.
   
   \eop
   
\bigskip

The decomposition of Proposition 
\ref{Prop: decomp no explosion} in terms of composition of linear functionals adapts easily also to the non-dynamical case: it holds for any matrix which admits a real logarithm. Further, for complex matrices the result holds with dimensions of the subspaces given by $d_i=1$ for all $i\in 1, \ldots, n$.

\end{document}